\input amstex
\documentstyle{amsppt}

\define\Cee{{\Bbb C}}

\define\Zee{\Bbb Z}

\define\Res#1{\operatorname{Res}_{#1 = 0}}
\define\hil#1{{X^{[#1]}}}

\define\W{\widetilde}
\define\w{\tilde}

\def\sumlim#1#2{\underset {#1} \to {\overset {#2} \to \sum}}

\def\bigopluslim#1#2{\underset {#1} \to {\overset {#2} \to \bigoplus}}

\define\End{\operatorname{End}}

\define\Id{\operatorname{Id}}

\define\Supp{\operatorname{Supp}}

\define\ch{\operatorname{ch}}
\define\td{\operatorname{td}}

\define\proof{\demo{Proof}}
\def\stop {\nobreak$\quad$\lower 1pt\vbox{
    \hrule
    \hbox to 7pt{\vrule height 7pt\hfil\vrule height 7pt}
      \hrule}\ifmmode\relax\else\par\medbreak \fi}
\define\endproof{\stop\enddemo}
\define\Remark#1{\medskip\noindent {\bf{Remark #1}.}}
\define\theorem#1{\medskip\noindent {\bf{Theorem #1}.} \it}
\define\lemma#1{\medskip\noindent {\bf{Lemma #1}.} \it}
\define\proposition#1{\medskip\noindent {\bf{Proposition #1}.} \it}
\define\corollary#1{\medskip\noindent {\bf{Corollary #1}.} \it}
\define\claim#1{\medskip\noindent {\bf{Claim #1}.} \it}
\define\endstatement{\rm }
%\define\theorem#1{\proclaim{Theorem #1}}
%\define\lemma#1{\proclaim{Lemma #1}}
%\define\proposition#1{\proclaim{Proposition #1}}
%\define\corollary#1{\proclaim{Corollary #1}}
%\define\claim#1{\proclaim{Claim #1}}

\define\section#1{\bigskip\noindent {\bf #1}}
%\define\section#1{\specialhead #1 \endspecialhead}
\define\ssection#1{\medskip\noindent{\bf #1}}

\def\today{\ifcase \month \or January\or February\or
  March\or April\or May\or June\or July\or August\or
  September\or October\or November\or December\fi\
  \number \day, \number \year}

\loadbold \hoffset=0.3 in \voffset=0.6 in

\leftheadtext{Wei-ping Li, Zhenbo Qin, and Weiqiang Wang}
\rightheadtext{Generators for the cohomology ring of Hilbert
schemes}

\topmatter
\title
Generators for the cohomology ring of Hilbert schemes of points on
surfaces
\endtitle
\author {Wei-ping Li, Zhenbo Qin, and Weiqiang Wang}
\endauthor
\address
Department of Mathematics, HKUST, Clear Water Bay, Kowloon, Hong
Kong
\endaddress
\email mawpli\@uxmail.ust.hk
\endemail
\address
Department of Mathematics, University of Missouri, Columbia, MO
65211, USA
\endaddress
\email zq\@math.missouri.edu
\endemail
\address
Department of Mathematics, NC State University, Raleigh, NC 27695;
Current address: Department of Mathematics, University of
Virginia, Charlottesville, VA 22904, USA \email ww9c\@virginia.edu
\endemail
\endaddress
%\keywords Hilbert schemes, projective surfaces, and cohomology
%rings
%\endkeywords
%\thanks
%Received 18 September 2000. Revision received 27 March 2001.
%\endthanks
%\subjclass Primary 14C05; Secondary 17B69
%\endsubjclass
\abstract Using the methods developed in \cite{LQW}, we obtain a
second set of generators for the cohomology ring of the Hilbert
scheme of points on an arbitrary smooth projective surface $X$
over the field of complex numbers. These generators have clear and
simple geometric as well as algebraic descriptions.
\endabstract
\endtopmatter

\NoBlackBoxes \TagsOnRight
\document
\section{1. Introduction}

In recent years, there has been a surge of research interest in
Hilbert schemes $\hil{n}$ of points on surfaces $X$ largely due to
the work of G\"ottsche \cite{Got}, Nakajima \cite{Na1}, Grojnowski
\cite{Gro}, and Vafa and Witten \cite{VW}. Earlier, Ellingsrud and
Str\o mme \cite{ES1} calculated the Betti numbers of the Hilbert
schemes of points on the projective plane $\Bbb P^2$, the affine
plane $\Bbb C^2$, and rational ruled surfaces. Subsequently,
G\"ottsche \cite{Got} determined the Betti numbers of $\hil{n}$
for any surfaces $X$ by finding a beautiful formula for the
generating function of the Betti numbers of $\hil{n}$.
G\"ottsche's results suggested that one should study Hilbert
schemes $\hil{n}$ for all $n$ altogether rather than study each
$\hil{n}$ individually. This idea is also echoed in Vafa and
Witten's work \cite{V-W} which stated that the generating function
of the Euler numbers of $\hil{n}$ is the partition function of
some physical theory. Later on, Nakajima \cite{Na1} constructed a
Heisenberg algebra action on the direct sum of cohomology rings of
$\hil{n}$ over all $n$. Similar results were obtained by
Grojnowski \cite{Gro}.
\medskip

The cohomology ring structure of Hilbert schemes $\hil{n}$ is a
more subtle issue than the Betti numbers, and it has also been
studied extensively. Many questions in enumerative geometry can be
interpreted as questions in cup product in the cohomology ring of
$\hil{n}$ (see \cite{EL}). Ellingsrud and Str\o mme in \cite{ES2}
studied the ring structure of the cohomology ring of $\hil{n}$
when $X$ is the projective plane $\Bbb P^2$ and the affine plane
$\Bbb C^2$. They found a set of generators of the cohomology ring
structure on $H^*(\hil{n})$ via the Chern classes of the
tautological rank-$n$ vector bundles coming from the universal
subscheme. This result was extended by Beauville \cite{Bea} to
other rational surfaces and ruled surfaces. Very recently, Markman
\cite{Mar} further extended the method to K3 surfaces. This
method, in the context of moduli spaces $\frak M$ of stable
sheaves, basically says that if one can express the class of the
diagonal in $\frak M\times\frak M$ in terms of the Chern classes
of a universal sheaf $\Cal E$ on $\frak M\times X$, then one can
show that the K\"unneth components of Chern classes of $\Cal E$
provide a set of ring generators. Fantechi and G\"ottsche in
\cite{FG} also found a set of generators for the cohomology ring
of $\hil{3}$ for any surface $X$. There is a totally different
approach initiated by Lehn \cite{Leh} where the Heisenberg algebra
construction mentioned earlier is used in an essential way. In
particular, Lehn was able to describe the cohomology ring of
$(\Bbb C^2)^{[n]}$ in terms of certain explicit differential
operators. There are also different viewpoints, as first indicated
in the work of I. Frenkel and the third author \cite{FW}, of
relating the cohomology rings of Hilbert schemes for the affine
plane $\Cee^2$ to the convolution product on symmetric groups (see
\cite{LS1, Vas}).

\medskip
In our paper \cite{LQW}, a set of $(n \cdot \dim H^*(X))$
generators for the cohomology ring of $\hil{n}$ for any smooth
projective surface $X$ was found. These generators, denoted by
$G_{i}(\gamma, n)$ where $0 \le i < n$ and $\gamma$ runs over a
linear basis of $H^*(X)$, were defined by using essentially the
Chern classes of the universal subscheme (see Definition 2.8 (vi)
for details). They coincide with the generators found by
Ellingsrud and Str\o mme in the case of $\Bbb P^2$. Note that
although the statement on this set of generators has nothing to do
with Heisenberg algebras etc, the machinery which we built up in
its proof is deeply rooted in the theory of vertex algebras and
the work of Nakajima, Grojnowski and Lehn. We took the viewpoint
effectively that the cup products with $G_{i}(\gamma, n)$ for all
$n$ associated to a fixed $\gamma$ should be treated as a single
operator acting on the direct sum of cohomology rings of $\hil{n}$
over all $n$.

\medskip
The geometric interpretation of the ring generators in terms of
the universal subscheme is not always clear however. In this
paper, we present a new set of $(n \cdot \dim H^*(X))$ generators
for the cohomology rings of Hilbert schemes, which has a simple
geometric interpretation. This new set of generators also affords
a nice algebraic interpretation in terms of the Heisenberg algebra
operators. To be precise, let $|0\rangle$ be the element $1$ of
$H^0(\hil{0})=\Bbb Q$. Denote by $\frak q_j(\gamma)$ the
Heisenberg algebra operators, where $j \in \Zee$ and $\gamma\in
H^*(X)$ (see Definition 2.8 (ii)). For $0 \le i < n$ and $\gamma
\in H^*(X)$, define a cohomology class $B_i(\gamma, n) \in
H^*(\hil{n})$ by putting $$B_i(\gamma, n) = {1 \over (n-i-1)!}
\cdot \frak q_{i+1}(\gamma) \frak q_{1}(1_X)^{n-i-1} |0\rangle.$$
Note that these are the simplest cohomology classes in
$H^*(\hil{n})$ both in geometric terms and Heisenberg algebraic
terms. Moreover, $B_0(1_X, n) = n \cdot 1_{\hil{n}}$. In addition,
if either $i > 0$ or $\gamma \in H^s(X)$ with $s > 0$, then one
can easily show that the Poincar\' e dual of $B_i(\gamma, n)$ is
the homology class represented by the closed subset: $$\{\, \xi\in
\hil{n}\, |\, \hbox{$\exists x\in \Gamma$ with $\ell(\xi_x)\ge
i+1$}\,\}$$ where $\Gamma$ is a homology cycle of $X$ representing
the Poincar\'e dual of $\gamma$, and $\xi_x$ is the component of
$\xi$ such that $\xi_x$ is supported at $x$.

The following is our main result in this paper.

\theorem{1.1} Let $X$ be a smooth projective surface over the
field of complex numbers. For $n \ge 1$, the cohomology ring $\Bbb
H_n = H^*(\hil{n})$ is generated by the classes $B_{i}(\gamma, n)$
where $0 \le i < n$ and $\gamma$ runs over a linear basis of
$H^*(X)$.
\endstatement

This theorem is proved in section 3 after we review some
definitions and results from \cite {Na2, Leh, LQW} in Section 2.
The main idea is that by exploring the results established in
\cite{LQW}, we are able to determine certain relations among the
cohomology classes $G_{i}(\gamma, n)$ introduced in \cite{LQW} and
$B_{i}(\gamma, n)$ introduced here. More precisely, we show that
$B_i(\gamma, n)$ coincides with $G_i(\gamma, n)$ for $i = 0, 1$.
In addition, for $2 \le i < n$, $B_i(\gamma, n)$ is equal to
$(-1)^i(i+1)! \cdot G_i(\gamma, n)$ plus a finite sum of products
of the form
$$G_{m_1}(\gamma_1, n) \cdot \, \ldots \, \cdot G_{m_t}(\gamma_t,
n) \in \Bbb H_n$$ where $m_1, \ldots, m_t \ge 0$ with $m_1 +
\ldots + m_t < i$, and $\gamma_1, \ldots, \gamma_t \in H^*(X)$. We
remark that in the process of deriving the relations among
$G_{i}(\gamma, n)$ and $B_{i}(\gamma, n)$, we actually find a new
proof to our theorem in \cite{LQW} which states that the classes
$G_{i}(\gamma, n)$, where $0 \le i < n$ and $\gamma$ runs over a
linear basis of $H^*(X)$, generate the cohomology ring of
$\hil{n}$. Then Theorem~1.1 is derived by using this theorem and
the relations between the cohomology classes $G_{i}(\gamma, n)$
and $B_{i}(\gamma, n)$.

\medskip\noindent
{\bf Conventions:} Throughout the paper, all cohomology rings are
in $\Bbb Q$-coefficients. The cup product between two cohomology
classes $\alpha$ and $\beta$ is denoted by $\alpha \cdot \beta$ or
simply by $\alpha \beta$. For a continuous map $p: Y_1 \to Y_2$
between two smooth compact manifolds and for $\alpha_1 \in
H^*(Y_1)$, the push-forward $p_*(\alpha_1)$ is defined by
$$p_*(\alpha_1) = \text{PD}^{-1}p_{*}(\text{PD}(\alpha_1))$$ where
$\text{PD}$ stands for the Poincar\'e duality. Unless otherwise
specified, we make no distinction between an algebraic cycle and
its corresponding cohomology class so that intersections among
algebraic cycles correspond to cup products among the
corresponding cohomology classes. For instance, for two algebraic
cycles $[a]$ and $[b]$ on a smooth projective variety $Y$, it is
understood that $[a] \cdot [b] \in H^*(Y)$.

\section{2. Results from \cite{Na2, Leh, LQW}}

In this section, we shall fix some notations, and recall some
results from \cite{Na2, Leh, LQW}. For convenience, we also review
certain basic facts for the Hilbert scheme of points in a smooth
projective surface.

Let $X$ be a smooth projective surface over $\Cee$, and $\hil{n}$
be the Hilbert scheme of points in $X$. An element in the Hilbert
scheme $\hil{n}$ is represented by a length-$n$ $0$-dimensional
closed subscheme $\xi$ of $X$, which sometimes is called a
length-$n$ $0$-cycle. For $\xi \in \hil{n}$, let $I_{\xi}$ and
$\Cal O_\xi$ be the corresponding sheaf of ideals and structure
sheaf respectively. For a point $x \in X$, let $\xi_x$ be the
component of $\xi$ supported at $x$ and $I_{\xi, x} \subset \Cal
O_{X, x}$ be the stalk of $I_{\xi}$ at $x$. It is known from
\cite{Fog} that $\hil{n}$ is smooth. In $\hil{n}\times X$, we have
the universal codimension-$2$ subscheme: $$\Cal Z_n=\{(\xi, x)
\subset \hil{n}\times X \, | \, x\in \Supp{(\xi)}\}\subset
\hil{n}\times X.  \tag 2.1$$ Also, we let
$X^n=\underbrace{X\times\cdots\times X}_{n \text{ times}}$ be the
$n$-th Cartesian product, and $$X^{[n_1], \ldots, [n_k]} =
\hil{n_1} \times \cdots \times \hil{n_k}. \tag 2.2$$

\medskip\noindent
{\bf Definition 2.3.} (i) Let $\Bbb H = \bigopluslim{n, i \ge 0}{}
\Bbb H^{n,i}$ denote the double-graded vector space with
components $\Bbb H^{n,i} \overset \text{def} \to = H^i(\hil{n})$,
and $\Bbb H_n \overset \text{def} \to = H^*(\hil{n}) \overset
\text{def} \to = \bigopluslim{i=0}{4n} H^i(\hil{n})$. The element
$1$ in $H^0(\hil{0}) = \Bbb Q$ is called the {\it vacuum vector}
and denoted by $|0\rangle$;

(ii) A linear operator $\frak f \in \End(\Bbb H)$ is {\it
homogeneous of bidegree} $(\ell, m)$ if $$\frak f(\Bbb H^{n,i})
\subset \Bbb H^{n+\ell,i+m}.  \tag 2.4$$ Furthermore, $\frak f \in
\End(\Bbb H)$ is {\it even} (resp., {\it odd}) if $m$ is even
(resp., odd).

(iii) For two homogeneous linear operators $\frak f$ and $\frak g
\in \End(\Bbb H)$ of bidegrees $(\ell, m)$ and $(\ell_1, m_1)$
respectively, define the {\it Lie superalgebra bracket} $[\frak f,
\frak g]$ by $$[\frak f, \frak g] = \frak f \frak g - (-1)^{m m_1}
\frak g \frak f.   \tag 2.5$$

(iv) Let $\frak A(\alpha), \frak B(\beta) \in \End(\Bbb H)$ be two
series of operators depending linearly on $\alpha, \beta \in
H^*(X)$. Then, the commutator $[\frak A(\alpha), \frak B(\beta)]$
satisfies {\it the transfer property} if $[\frak A(\alpha), \frak
B(\beta)] = [\frak A(1_X), \frak B(\alpha \beta)] = [\frak
A(\alpha \beta), \frak B(1_X)]$ for all $\alpha, \beta$.

A non-degenerate super-symmetric bilinear form $(, )$ on $\Bbb H$
is induced from the standard one on $\Bbb H_n=H^*(\hil{n})$. For a
homogeneous linear operator $\frak f \in \End(\Bbb H)$ of bidegree
$(\ell, m)$, we can define its {\it adjoint} $\frak f^\dagger \in
\End(\Bbb H)$ by $$(\frak f(\alpha), \beta) = (-1)^{m \cdot
|\alpha|} \cdot (\alpha, \frak f^\dagger(\beta)) \tag 2.6$$ where
$|\alpha| = s$ for $\alpha \in H^s(X)$. Note that the bidegree of
$\frak f^\dagger$ is $(-\ell, m - 4 \ell)$. Also,
$$(\frak f \frak g)^\dagger = (-1)^{m m_1} \cdot \frak g^\dagger
\frak f^\dagger \qquad \text{and} \qquad [\frak f, \frak
g]^\dagger = -[\frak f^\dagger, \frak g^\dagger] \tag 2.7$$ where
$\frak g \in \End(\Bbb H)$ is another homogeneous linear operator
of bidegree $(\ell_1, m_1)$.

Next, we collect from \cite{Na2, Leh, LQW} the definitions of the
closed subset $Q^{[n+\ell,n]}$ in $\hil{n+\ell} \times X \times
\hil{n}$, the Heisenberg generator $\frak q_n$, the Virasoro
generator $\frak L_n$, the boundary operator $\frak d$, the
derivative $\frak f'$ of a linear operator $\frak f \in \End(\Bbb
H)$, and the operators $\frak G_i(\alpha) \in \End(\Bbb H)$ for $i
\ge 0$ and $\alpha \in H^*(X)$.

\medskip\noindent
{\bf Definition 2.8.} (i) For $n \ge 0$, define $Q^{[n,n]} =
\emptyset$. For $n \ge 0$ and $\ell > 0$, define $Q^{[n+\ell,n]}
\subset \hil{n+\ell} \times X \times \hil{n}$ to be the closed
subset $$\{ (\xi, x, \eta) \in \hil{n+\ell} \times X \times
\hil{n} \, | \, \xi \supset \eta \text{ and } \Supp(I_\eta/I_\xi)
= \{ x \} \}; \tag 2.9$$

(ii) For $n \in \Zee$, define linear maps $\frak q_n: H^*(X) \to
\End(\Bbb H)$ as follows. When $n \ge 0$, the linear operator
$\frak q_n(\alpha) \in \End(\Bbb H)$ with $\alpha \in H^*(X)$ is
defined by $$\frak q_n(\alpha)(a) = \w p_{1*}([Q^{[m+n,m]}] \cdot
\w \rho^*\alpha \cdot \w p_2^*a)  \tag 2.10$$ for all $a \in \Bbb
H_m = H^*(\hil{m})$, where $[Q^{[m+n,m]}]$ is (the cohomology
class corresponding to) the algebraic cycle associated to
$Q^{[m+n,m]}$, and $\w p_1, \w \rho, \w p_2$ are the projections
of $\hil{m+n} \times X \times \hil{m}$ to $\hil{m+n}, X, \hil{m}$
respectively. When $n < 0$, define the operator $\frak q_n(\alpha)
\in \End(\Bbb H)$ with $\alpha \in H^*(X)$ by $$\frak q_n(\alpha)
= (-1)^n \cdot \frak q_{-n}(\alpha)^\dagger; \tag 2.11$$

(iii) For $n \in \Zee$, define linear maps $\frak L_n: H^*(X) \to
\End(\Bbb H)$ by putting $$\frak L_n = \cases
  {1 \over 2} \cdot \sumlim{m \in \Zee}{} \frak q_m \frak q_{n-m}
     \tau_{2*}, &\text{if $n \ne 0$}  \\
  \sumlim{m > 0}{} \frak q_m \frak q_{-m} \tau_{2*}, &\text{if $n = 0$}
\\
\endcases
\tag 2.12$$ where $\tau_{2*}: H^*(X) \to H^*(X^2)$ is the linear
map induced by the diagonal embedding $\tau_2: X \to X^2$, and the
operator $\frak q_m \frak q_{\ell} \tau_{2*}(\alpha)$ stands for
$$\sum_j \frak q_m(\alpha_{j,1}) \frak q_{\ell}(\alpha_{j,2}) \tag
2.13$$ when $\tau_{k*}\alpha = \sum_j \alpha_{j,1} \otimes
\alpha_{j, 2}$ via the K\"unneth decomposition of $H^*(X^2)$;

(iv) Define the linear operator $\frak d \in \End(\Bbb H)$ by
$$\frak d = \bigopluslim{n}{} c_1(p_{1*}\Cal O_{\Cal Z_n}) =
\bigopluslim{n}{} (-[\partial \hil{n}]/2)   \tag 2.14$$ where
$p_1$ is the projection of $\hil{n} \times X$ to $\hil{n}$,
$\partial \hil{n}$ is the boundary of $\hil{n}$ consisting of all
$\xi \in \hil{n}$ with $|\Supp(\xi)| < n$, and the first Chern
class $c_1(p_{1*}\Cal O_{\Cal Z_n})$ of the rank-$n$ bundle
$p_{1*}\Cal O_{\Cal Z_n}$ acts on $\Bbb H_n = H^*(\hil{n})$ by the
cup product.

(v) For a linear operator $\frak f \in \End(\Bbb H)$, define its
{\it derivative} $\frak f'$ by $$\frak f' \overset \text{def} \to
= [\frak d, \frak f].  \tag 2.15$$ The higher derivative $\frak
f^{(k)}$ of $\frak f$ is defined inductively by $\frak f^{(k)} =
[\frak d, \frak f^{(k-1)}]$.

(vi) For $\alpha \in H^*(X)$ and $n \ge 0$, let $G_i(\alpha, n)$
be the $H^{|\alpha|+2i}(\hil{n})$-component of $$p_{1*}(\ch(\Cal
O_{\Cal Z_n}) \cdot p_2^*\td(X) \cdot p_2^*\alpha) \in \Bbb H_n$$
where $p_2$ is the projection of $\hil{n} \times X$ to $X$. For $i
\ge 0$ and $\alpha \in H^*(X)$, define $\frak G_i(\alpha) \in
\End({\Bbb H})$ to be the operator which acts on the component
$\Bbb H_n = H^*(\hil{n})$ by the cup product by the cohomology
class $G_i(\alpha, n)$.

\theorem{2.16} Let $K_X$ and $c_2(X)$ be the canonical divisor and
the second Chern class of $X$ respectively. Let $k \ge 0, n,m \in
\Zee$ and $\alpha, \beta \in H^*(X)$. Then,
\medskip
{\rm (i)} $[\frak q_n(\alpha), \frak q_m(\beta)] = n \cdot
\delta_{n+m} \cdot \int_X(\alpha \beta) \cdot \Id_{\Bbb H}$;
\medskip
{\rm (ii)} $[\frak L_n(\alpha), \frak q_m(\beta)] = -m \cdot \frak
q_{n+m}(\alpha \beta)$;
\medskip
{\rm (iii)} $[\frak L_n(\alpha), \frak L_m(\beta)] = (n-m) \cdot
\frak L_{n+m}(\alpha \beta) - {n^3-n \over 12} \cdot \delta_{n+m}
\cdot \int_X(c_2(X) \alpha \beta) \cdot \Id_{\Bbb H}$;
\medskip
{\rm (iv)} $\frak q_n'(\alpha) = n \cdot \frak L_n(\alpha) +
{n(|n|-1) \over 2} \frak q_n(K_X \alpha)$;
\medskip
{\rm (v)} $[\frak G_k(\alpha), \frak q_1(\beta)] = {1 \over k!}
\cdot \frak q_1^{(k)}(\alpha \beta)$;
\medskip
{\rm (vi)} $[\ldots [\frak G_k(\alpha), \frak q_{n_1}(\alpha_1)],
\ldots], \frak q_{n_{k+1}}(\alpha_{k+1})]$ is equal to
$$-\prod_{\ell = 1}^{k+1} (-n_\ell) \cdot \frak q_{n_1 + \ldots
+n_{k+1}}(\alpha \alpha_1 \cdots \alpha_{k+1})$$ for all $n_1,
\ldots , n_{k+1} \in \Zee$ with $\sum_{\ell =1}^{k+1} n_\ell \ne
0$ and all $\alpha_1, \ldots, \alpha_{k+1} \in H^*(X)$.
\endstatement

We notice that Theorem 2.16 (i) was proved by Nakajima \cite{Na1}
subject to some universal nonzero constant, which was determined
subsequently in [ES3]. The next three formulas in Theorem 2.16
were obtained by Lehn \cite{Leh}. Theorem 2.16 (v) is the Lemma
5.8 in \cite{LQW}, which is a generalization of a remarkable
theorem of Lehn (\cite{Leh}, Theorem~4.2). Theorem 2.16 (vi)
follows from the Theorem 5.13 (ii) and Proposition 4.10 (i) in
\cite{LQW}, and the proof of it uses the part (v) above. Also, as
observed by Nakajima and Grojnowski in \cite{Na1, Gro}, $\Bbb H$
is an irreducible representation of the Heisenberg algebra
generated by the $\frak q_i(\alpha)$'s with $|0\rangle \in
H^0(X^{[0]})$ being the highest weight vector. So a linear basis
of $\Bbb H$ is given by $\frak q_{i_1}(\alpha_1) \frak
q_{i_2}(\alpha_2) \cdots \frak q_{i_k}(\alpha_k) |0\rangle$ where
$k \ge 0$, $i_1 \geq i_2 \geq \cdots \geq i_k >0$, and each of
$\alpha_1, \alpha_2, \ldots, \alpha_k$ runs over a fixed linear
basis of $H^*(X)$.

The following elementary but convenient result is the Lemma 5.26
in \cite{LQW}.

\lemma{2.17} Fix $a, b$ with $1 \le a \le b$. Let $\frak g \in
\End(\Bbb H)$ be of bidegree $(\ell, s)$, and $$A = \frak
q_{m_1}(\beta_{1}) \cdots \frak q_{m_b}(\beta_{b}) |0\rangle.$$
Then, $\frak g(A)$ is equal to the sum of the following two terms:
$$\sum_{i=0}^{a-1} \sum_{\sigma_i} \pm
 \prod_{\ell \in \sigma_i^0} \frak q_{m_\ell}(\beta_{\ell})
 [\cdots [\frak g, \frak q_{m_{\sigma_i(1)}}(\beta_{\sigma_i(1)})],
 \cdots], \frak q_{m_{\sigma_i(i)}}(\beta_{\sigma_i(i)})]
 |0\rangle  \tag 2.18$$
and $$\sum_{\sigma_a} (-1)^{\sum_{k=0}^{a-1}
 (s+\sum_{\ell=1}^k |\beta_{\sigma_a(\ell)}|)
 \sum_{\sigma_a(k)<j<\sigma_a(k+1)}|\beta_j|}
 \prod_{\ell \in \sigma_a^1} \frak q_{m_\ell}(\beta_{\ell}) \cdot$$
$$\cdot  [\cdots [\frak g,
 \frak q_{m_{\sigma_a(1)}}(\beta_{\sigma_a(1)})],
 \cdots], \frak q_{m_{\sigma_a(a)}}(\beta_{\sigma_a(a)})]
 \prod_{\ell \in \sigma_a^2} \frak q_{m_\ell}(\beta_{\ell})
 |0\rangle  \tag 2.19$$
where for each fixed $i$ with $0 \le i \le a$, $\sigma_i$ runs
over all the maps $$\{\, 1, \ldots, i \,\} \to \{\, 1, \ldots, b
\,\}$$ satisfying $\sigma_i(1) < \cdots < \sigma_i(i)$. Moreover,
$\sigma_i^0 = \{ \ell \, | \, 1 \le \ell \le b, \ell \ne
\sigma_i(1), \ldots, \sigma_i(i) \}$, $\sigma_a^1 = \{ \ell \, |
\, 1 \le \ell < \sigma_a(a), \ell \ne \sigma_a(1), \ldots,
\sigma_a(a) \}$, and $\sigma_a^2 = \{ \ell \, | \, \sigma_a(a) <
\ell \le b \}$.
\endstatement
\proof Follows from moving all the commutators $$[\cdots [\frak g,
\frak q_{m_{\sigma_i(1)}}(\beta_{\sigma_i(1)})], \cdots], \frak
q_{m_{\sigma_i(i)}}(\beta_{\sigma_i(i)})]$$ with $0 \le i \le
(a-1)$ to the right, and applying the fact that $$\frak g_1 \frak
g_2 = [\frak g_1, \frak g_2] + (-1)^{m_1 m_2} \frak g_2 \frak
g_1$$ for $\frak g_1, \frak g_2 \in \End(\Bbb H)$ of bidegrees
$(\ell_1, m_1), (\ell_2, m_2)$ respectively.
\endproof

\section{3. A second set of generators for the cohomology ring}

In this section, we shall prove that the $(n \cdot \dim H^*(X))$
cohomology classes $B_i(\gamma, n)$ in Theorem 1.1 generate the
cohomology ring $\Bbb H_n = H^*(\hil{n})$. Moreover, we obtain an
alternative proof to the Theorem 5.30 in \cite{LQW} which has been
proved there by a different method. Finally, letting $\frak
B_i(\gamma) \in \End(\Bbb H)$ be the operator defined in (3.3), we
show that the commutator $[\frak B_i(\gamma), \frak q_n(\alpha)]$
satisfies the transfer property (see Definition 2.3 (iv) for the
definition of the transfer property).

\medskip\noindent
{\bf Definition 3.1.} (i) For $0 \le i < n$ and $\gamma \in
H^*(X)$, define $$B_i(\gamma, n) = {1 \over (n-i-1)!} \cdot \frak
q_{i+1}(\gamma) \frak q_{1}(1_X)^{n-i-1} |0\rangle \in \Bbb H_n;
\tag 3.2$$

(ii) For $i \ge 0$ and $\gamma \in H^*(X)$, define the operator
$\frak B_i(\gamma) \in \End(\Bbb H)$ by $$\frak B_i(\gamma) =
\bigopluslim{n \ge 0}{} B_i(\gamma, n) \tag 3.3$$ where
$B_i(\gamma, n)$ acts on the component $\Bbb H_n = H^*(\hil{n})$
by the cup product.

Notice that $\frak B_i(\gamma) \in \End(\Bbb H)$ is homogeneous of
bidegree $(0, |\gamma|+2i)$. Also, $$\frak B_i(\gamma)' = 0 \qquad
\text{and} \qquad \frak B_i(\gamma)^\dagger = \frak B_i(\gamma).
\tag 3.4$$ Our goal is to show that the cohomology ring $\Bbb H_n
= H^*(\hil{n})$ is generated by the $(n \cdot \dim H^*(X))$
classes $B_{i}(\gamma, n)$ where $0 \le i < n$ and $\gamma$ runs
over a linear basis of $H^*(X)$. We shall start with a technical
lemma which allows us to present an alternative proof to the
Theorem 5.30 in \cite{LQW} (see Theorem 3.19 below). Recall the
cohomology class $G_{i}(\gamma, n)$ defined in Definition 2.8
(vi).

\lemma{3.5} For $1 \le i \le n$, let $\Bbb H_{n, i}'$ be the
linear span of all the classes $$G_{m_1}(\gamma_1, n) \cdot \,
\ldots \, \cdot G_{m_t}(\gamma_t, n) \in \Bbb H_n   \tag 3.6$$
where $m_1, \ldots, m_t \ge 0$ with $m_1 + \ldots + m_t < i$, and
$\gamma_1, \ldots, \gamma_t \in H^*(X)$. Then, $$\frak
q_{1}(1_X)^{n-i} \frak q_{n_1}(\alpha_1) \cdots \frak
q_{n_k}(\alpha_k) |0\rangle \in \Bbb H_{n, i}'  \tag 3.7$$ for all
positive integers $n_1, \ldots, n_k$ with $\sumlim{\ell=1}{k}
n_\ell = i$ and all $\alpha_1, \ldots, \alpha_k \in H^*(X)$.
\proof We shall use induction on $i$. First of all, assume $i=1$.
Then, $k=1$ and $n_1 = 1$. By Theorem 2.16 (v), $[\frak
G_0(\alpha_1), \frak q_1(1_X)] = \frak q_1(\alpha_1)$. Thus, we
have $$\align &\qquad (n-1)! \cdot G_0(\alpha_1, n)
 = (n-1)! \cdot \frak G_0(\alpha_1)1_{\hil{n}}       \\
&= {1 \over n} \cdot \frak G_0(\alpha_1) \frak q_1(1_X)^n
|0\rangle
 = \frak q_{1}(1_X)^{n-1} \frak q_1(\alpha_1) |0\rangle \tag 3.8
\endalign$$
noting that $1_{\hil{n}} = {1 \over n!} \cdot \frak q_1(1_X)^n
|0\rangle$ and $\frak G_0(\alpha_1)|0\rangle = 0$. So (3.7) is
true when $i=1$.

Next, fixing an integer $i_0$ satisfying $1 \le i_0 < n$, we
assume that (3.7) holds for all the integers $i$ with $1 \le i \le
i_0$. We shall prove that (3.7) holds as well for $i = (i_0+1)$.
In other words, we shall verify that $$\frak
q_{1}(1_X)^{n-(i_0+1)} \frak q_{n_1}(\alpha_1) \cdots \frak
q_{n_k}(\alpha_k) |0\rangle \in \Bbb H_{n, i_0+1}'  \tag 3.9$$ for
all positive integers $n_1, \ldots, n_k$ with $\sumlim{\ell=1}{k}
n_\ell = i_0+1$ and all $\alpha_1, \ldots, \alpha_k \in H^*(X)$.

Let us explain how our induction works before we get into the
details. We shall take a suitably chosen element $A$ which lies in
$\Bbb H'_{n,i_0}$, and consider $\frak G_m(\alpha)(A) =
G_m(\alpha, n) \cdot A$ which lies in $\Bbb H'_{n,i_0 +1}$ for
some suitably chosen $m \ge 0$ and $\alpha \in H^*(X)$. By
applying Lemma 2.17 to the situation at hand, we shall observe
that the summation (2.18) in $\frak G_m(\alpha)(A)$ is in $\Bbb
H'_{n,i_0}$, and all the terms in (2.19) except those which
coincide with (3.9) are also in $\Bbb H'_{n,i_0}$. Since $\Bbb
H'_{n,i_0} \subset \Bbb H'_{n,i_0+1}$, (3.9) follows.

Now let us be more precise. Since $\Bbb H_{n, i_0}' \subset \Bbb
H_{n, i_0+1}'$, by induction hypothesis, we may assume that $\frak
q_{n_\ell}(\alpha_\ell) \ne \frak q_{1}(1_X)$ for all the integers
$\ell$ with $1 \le \ell \le k$. Put $$\align
 A
&=\frak q_{1}(1_X)^{n-(i_0+1-n_1)} \frak q_{n_2}(\alpha_2) \cdots
  \frak q_{n_k}(\alpha_k) |0\rangle  \\
&{\overset \text{def} \to =}
  \frak q_{m_1}(\beta_1) \cdots \frak q_{m_b}(\beta_b) |0\rangle
  \tag 3.10
\endalign$$
where we have $b = (n-(i_0+1-n_1)) + (k-1)$, and $$\align &\frak
q_{m_\ell}(\beta_\ell) = \frak q_{1}(1_X) \qquad
  \text{for } 1 \le \ell \le n-(i_0+1-n_1)  \tag 3.11   \\
&\frak q_{m_\ell}(\beta_\ell)
  \ne \frak q_{1}(1_X) \qquad
  \text{for } n-(i_0+1-n_1) < \ell \le b .  \tag 3.12
\endalign$$

Consider $G_{n_1-1}(\alpha_1, n) \cdot A$. Note that $0 \le
(i_0+1-n_1) \le i_0$. If $(i_0+1-n_1)=0$, then $A = \frak
q_{1}(1_X)^n|0\rangle = (n-1)! \cdot G_0(1_X, n) \in \Bbb H_{n,
1}'$ by (3.8); so $$G_{n_1-1}(\alpha_1, n) \cdot A \in
G_{n_1-1}(\alpha_1, n) \cdot \Bbb H_{n, 1}' \subset \Bbb H_{n,
n_1}' = \Bbb H_{n, i_0+1}'.  \tag 3.13$$ If $(i_0+1-n_1) > 0$,
then by induction hypothesis, $A \in \Bbb H_{n, i_0+1-n_1}'$; so
$$G_{n_1-1}(\alpha_1, n) \cdot A \in G_{n_1-1}(\alpha_1, n) \cdot
\Bbb H_{n, i_0+1-n_1}' \subset \Bbb H_{n, i_0}' \subset \Bbb H_{n,
i_0+1}'. \tag 3.14$$ In summary, we have showed that
$G_{n_1-1}(\alpha_1, n) \cdot A \in \Bbb H_{n, i_0+1}'$. Thus,
$$\frak G_{n_1-1}(\alpha_1)(A) =G_{n_1-1}(\alpha_1, n) \cdot A \in
\Bbb H_{n, i_0+1}'.  \tag 3.15$$

Applying Lemma 2.17 to $a = n_1$ and $\frak g = \frak
G_{n_1-1}(\alpha_1) = \frak G_{a-1}(\alpha_1)$, we see that the
class $\frak G_{a-1}(\alpha_1)(A)$ consists of two parts (2.18)
and (2.19). By (3.11), the number of $\frak q_{1}(1_X)$'s in every
nonvanishing term of (2.18) is at least $$(n-(i_0+1-n_1)) - (a -1)
= n - i_0.  \tag 3.16$$ So by induction hypothesis, (2.18) is
contained in $\Bbb H_{n, i_0}'$. Let $N(\sigma_a)$ be the number
of $\frak q_{1}(1_X)$'s in a nonvanishing term in (2.19)
corresponding to $\sigma_a$. By (3.11) again, $$N(\sigma_a) \ge
(n-(i_0+1-n_1)) - a = n - (i_0+1).  \tag 3.17$$ If $N(\sigma_a) >
n - (i_0+1)$, then by induction, this nonvanishing term in (2.19)
corresponding to $\sigma_a$ is contained in $\Bbb H_{n, i_0}'$.
Moreover, we see from (3.12) that $$N(\sigma_a) = n - (i_0+1)$$ if
and only if $\sigma_a(1) < \ldots < \sigma_a(a) \le
n-(i_0+1-n_1)$. So this nonvanishing term $$(-1)^{\sum_{k=0}^{a-1}
 (s+\sum_{\ell=1}^k |\beta_{\sigma_a(\ell)}|)
 \sum_{\sigma_a(k)<j<\sigma_a(k+1)}|\beta_j|}
 \prod_{\ell \in \sigma_a^1} \frak q_{m_\ell}(\beta_{\ell}) \cdot$$
$$\cdot  [\cdots [\frak G_{a-1}(\alpha_1),
 \frak q_{m_{\sigma_a(1)}}(\beta_{\sigma_a(1)})],
 \cdots], \frak q_{m_{\sigma_a(a)}}(\beta_{\sigma_a(a)})]
 \prod_{\ell \in \sigma_a^2} \frak q_{m_\ell}(\beta_{\ell})
 |0\rangle$$
in (2.19) corresponding to $\sigma_a$ can be simplified to $$\frak
q_1(1_X)^{|\sigma_a^1|} [\cdots [\frak G_{a-1}(\alpha_1),
\underbrace{\frak q_1(1_X)], \cdots], \frak q_1(1_X)]}_{a \text{
times}} \prod_{\ell \in \sigma_a^2} \frak q_{m_\ell}(\beta_{\ell})
|0\rangle.  $$ By Theorem 2.16 (vi), the above term can be further
simplified to $$\align &\qquad \frak q_1(1_X)^{|\sigma_a^1|}
(-1)^{a+1} \frak q_a(\alpha_1)
 \prod_{\ell \in \sigma_a^2} \frak q_{m_\ell}(\beta_{\ell}) |0\rangle
\\
&= (-1)^{n_1+1} \cdot
 \frak q_{1}(1_X)^{n-(i_0+1)} \frak q_{n_1}(\alpha_1) \cdots
 \frak q_{n_k}(\alpha_k) |0\rangle.
\endalign$$
Now there are ${n-(i_0+1-n_1) \choose a} = {n-(i_0+1-n_1) \choose
n_1}$ such terms in (2.19). Therefore, $$\align &\quad \frak
G_{n_1-1}(\alpha_1)(A) = \frak G_{a-1}(\alpha_1)(A)  \\ &\equiv
{n-(i_0+1)+n_1 \choose n_1} \cdot (-1)^{n_1+1} \cdot  \\
&\qquad\qquad \cdot \frak q_{1}(1_X)^{n-(i_0+1)} \frak
q_{n_1}(\alpha_1) \cdots \frak q_{n_k}(\alpha_k) |0\rangle   \tag
3.18
\endalign$$
modulo ${\Bbb H_{n, i_0}'}$. Since $\Bbb H_{n, i_0}' \subset \Bbb
H_{n, i_0+1}'$, (3.9) follows from (3.15) and (3.18).
\endproof

\theorem{3.19} For $n \ge 1$, the cohomology ring $\Bbb H_n =
H^*(\hil{n})$ is generated by the classes $G_{i}(\gamma, n)$ where
$0 \le i < n$ and $\gamma$ runs over a linear basis of $H^*(X)$.
\endstatement
\proof Follows immediately from (3.7) by taking $i = n$. We remark
that this is the Theorem 5.30 in \cite{LQW}, and has been proved
there by a different method.
\endproof

The next lemma provides a relation between the classes
$B_i(\gamma, n)$ and $G_i(\gamma, n)$.

\lemma{3.20} {\rm (i)} For $\gamma \in H^*(X)$, we have
$$B_0(\gamma, n) = G_0(\gamma, n) \quad \text{and} \quad
B_1(\gamma, n) = -2G_1(\gamma, n);$$

{\rm (ii)} For $2 \le i < n$ and $\gamma \in H^*(X)$, we have
$$B_i(\gamma, n) \equiv (-1)^i(i+1)! \cdot G_i(\gamma, n) \pmod
{\Bbb H_{n, i}'}.$$
\endstatement
\proof (i) Recall from (3.2) that for $0 \le i < n$, we have
$$(n-i-1)! \cdot B_i(\gamma, n) = \frak q_{i+1}(\gamma) \frak
q_{1}(1_X)^{n-i-1} |0\rangle.  \tag 3.21$$ Combining (3.21) and
(3.8), we conclude that $$B_0(\gamma, n) = {1 \over (n-1)!} \cdot
\frak q_{1}(\gamma) \frak q_{1}(1_X)^{n-1} |0\rangle = G_0(\gamma,
n).$$

Next, we apply Lemma 2.17 to $\frak g = \frak G_1(\gamma)$, $A =
\frak q_{1}(1_X)^n |0\rangle$, $b=n$, and $a =2$. So $\frak
G_1(\gamma)(A)$ consists of two parts (2.18) and (2.19). By
Theorem 2.16 (v), $$[\frak G_1(\gamma), \frak q_{1}(1_X)] = \frak
q_{1}'(\gamma) = \frak L_1(\gamma).$$ Since $\frak
G_1(\gamma)|0\rangle =0$ and $\frak L_1(\gamma)|0\rangle =0$,
(2.18) is zero. By Theorem 2.16 (ii), $[[\frak G_1(\gamma), \frak
q_{1}(1_X)], \frak q_{1}(1_X)] =[\frak L_1(\gamma), \frak
q_{1}(1_X)] =- \frak q_2(\gamma)$. So (2.19) is equal to $$-{n
\choose 2} \cdot \frak q_2(\gamma) \frak q_{1}(1_X)^{n-2}
|0\rangle = -{n! \over 2} \cdot B_1(\gamma, n)$$ where we have
used (3.21). Therefore, we obtain $\frak G_1(\gamma)(A) = -{n!/2}
\cdot B_1(\gamma, n)$. Since $\frak G_1(\gamma)(A) = G_1(\gamma,
n) \cdot A = n! \cdot G_1(\gamma, n)$, we see that $B_1(\gamma, n)
= -2G_1(\gamma, n)$.

(ii) In the proof of Lemma 3.5, we take $\alpha_1 = \gamma$, $i_0
= i$, $a = n_1 = i+1$, $k=1$, $A = \frak q_{1}(1_X)^n |0\rangle$,
and $b=n$. We see from (3.18) that $$\frak G_i(\gamma)(A) \equiv
{n \choose i+1} \cdot (-1)^i \cdot \frak q_{1}(1_X)^{n-i-1} \frak
q_{i+1}(\gamma) |0\rangle \pmod {\Bbb H_{n, i}'}.  \tag 3.22$$
Since $n! \cdot 1_{\hil{n}} = A$, we conclude from (3.22) and
(3.21) that $$\align
 n! \cdot G_i(\gamma, n)
&= G_i(\gamma, n) \cdot A = \frak G_i(\gamma)(A)      \\ &\equiv
{n \choose i+1} \cdot (-1)^i \cdot
 \frak q_{1}(1_X)^{n-i-1} \frak q_{i+1}(\gamma) |0\rangle
 \pmod {\Bbb H_{n, i}'}    \\
&\equiv {n! \over (i+1)!} \cdot (-1)^i \cdot B_i(\gamma, n)
 \pmod {\Bbb H_{n, i}'}.
\endalign$$
It follows that $B_i(\gamma, n) \equiv (-1)^i(i+1)! \cdot
G_i(\gamma, n) \pmod {\Bbb H_{n, i}'}$.
\endproof

\theorem{3.23} For $n \ge 1$, the cohomology ring $\Bbb H_n =
H^*(\hil{n})$ is generated by $$B_{i}(\gamma, n) = \frak
B_{i}(\gamma)(1_{\hil{n}}) \tag 3.24$$ where $0 \le i < n$ and
$\gamma$ runs over a linear basis of $H^*(X)$. Moreover, the
relations among these generators are precisely the relations among
the restrictions $\frak B_i(\gamma)|_{\Bbb H_n}$ of the
corresponding operators $\frak B_i(\gamma)$ to ${\Bbb H_n}$.
\endstatement
\proof Note that the second statement follows from the fact that
the operators $\frak B_i(\gamma)|_{\Bbb H_n}$ are defined in terms
of the cup products by the cohomology classes $B_i(\gamma, n)$. In
the following, we prove the first statement.

Let $\W {\Bbb H}_n$ be the subring generated by the $(n \cdot \dim
H^*(X))$ cohomology classes in (3.24). In view of Theorem 3.19, it
suffices to show that $$G_{i}(\gamma, n) \in \W {\Bbb H}_n  \tag
3.25$$ for all $i$ with $0 \le i < n$ and all $\gamma \in H^*(X)$.
We shall use induction on $i$. First of all, this is true for $i =
0, 1$ by Lemma 3.20 (i). Next, fixing an integer $i_0$ with $2 \le
i_0 < n$, we assume that (3.25) is true for all the integers $i$
with $0 \le i < i_0$. We want to show that (3.25) holds as well
for $i = i_0$, i.e., $$G_{i_0}(\gamma, n) \in \W {\Bbb H}_n$$ for
all $\gamma \in H^*(X)$. Indeed, we see from Lemma 3.20 (ii) that
$$G_{i_0}(\gamma, n) \equiv {(-1)^{i_0} \over ({i_0}+1)!} \cdot
B_{i_0}(\gamma, n)  \pmod {\Bbb H_{n, {i_0}}'}.  \tag 3.26$$ By
the definition of $\Bbb H_{n, {i_0}}'$ in Lemma 3.5 and by the
induction hypothesis, we have $\Bbb H_{n, {i_0}}' \subset \W {\Bbb
H}_n$. So it follows from (3.26) that $G_{i_0}(\gamma, n) \in \W
{\Bbb H}_n$.
\endproof

We stress that Theorem~3.23 is not a mere consequence of
Theorem~3.19 which was first established in \cite{LQW}. As we
worked on the only way that we found to derive Theorem~3.23, we
happened to obtain a new proof of the old Theorem~3.19 which we
present here.

In the last part of the paper, we shall establish the {\it
transfer property} for the commutator $[\frak B_i(\gamma), \frak
q_n(\alpha)]$, cf. Definition 2.3 (iv). This property tells us how
for a fixed $i \ge 0$, the operators $\frak B_i(\gamma)$
associated with different classes $\gamma \in H^*(X)$ are related
to each other. The transfer property, which seems to be universal
for these types of operators in Hilbert schemes, was formulated
and emphasized in \cite{LQW}, and examples of such property first
appeared in \cite{Leh}. Such a property can be used as a tool to
show that various statements in Hilbert schemes are insensitive to
the underlying surface $X$ and can in principle be reduced to the
understanding of the affine plane case. We mention that the
transfer property for the commutator $[\frak G_i(\gamma), \frak
q_n(\alpha)]$ established in \cite{LQW} has been used effectively
in a new remarkable paper of Lehn and Sorger \cite{LS2}. We expect
that the transfer property for $[\frak B_i(\gamma), \frak
q_n(\alpha)]$ which we prove here will play an important role in a
further development.

We start with some notation, and prove a technical lemma which
gives an alternative description of the cohomology class
$B_i(\gamma, n)$. For $0 \le i < n$, define $$Z_{n, i+1} = \{(\xi,
x) \in \hil{n} \times X \, | \, \ell(\xi_x) \ge i + 1 \}   \tag
3.27$$ where $\xi_x$ is the component of $\xi$ such that $\xi_x$
is supported at $x$. Note that $Z_{n, i+1}$ is closed,
irreducible, and of dimension $(2n-i)$.

\lemma{3.28} For $0 \le i < n$ and $\gamma \in H^*(X)$, we have
$$B_i(\gamma, n) = p_{1*}([\Cal Z_{n, i+1}] \cdot p_2^*\gamma)$$
where $p_1$ and $p_2$ are the projections of $\hil{n} \times X$ to
$\hil{n}$ and $X$ respectively.
\endstatement
\proof First of all, let $i=0$ and $\gamma = 1_X$. By (3.2), we
have $$B_0(1_X, n) = {1 \over (n-1)!} \cdot \frak
q_{1}(1_X)^n|0\rangle = n \cdot 1_{\hil{n}} = p_{1*}([\Cal Z_{n}])
= p_{1*}([\Cal Z_{n, 1}])$$ noting that $1_{\hil{n}} = 1/n! \cdot
\frak q_{1}(1_X)^n|0\rangle$. So Lemma 3.28 holds for $i=0$ and
$\gamma = 1_X$. In fact, by linearity, Lemma 3.28 is true for $i =
0$ and $\gamma \in H^0(X)$.

Next, we assume that either $i > 0$ or $\gamma \in H^s(X)$ with $s
> 0$. Under these conditions, the Poincar\'e duals of the
cohomology classes $B_i(\gamma, n)$ and $$p_{1*}([\Cal Z_{n, i+1}]
\cdot p_2^*\gamma)$$ have clear geometric interpretations. So we
shall work in the homology setting (only in this proof) with the
help of Poincar\'e duality. Notice that via Poincar\'e duality,
the cup product on cohomology theory become the intersection in
homology theory, and the map $p_{1*}$ we defined in the
Conventions for cohomology theory becomes the ordinary pushforward
map $p_{1*}$ for homology theory. Let $\Gamma$ be a
$(4-s)$-dimensional homology cycle of $X$ representing the
Poincar\'e dual of $\gamma \in H^s(X)$, and $X^{[m]}_0$ be the
open dense subset of $\hil{m}$ consisting of all $\xi \in \hil{m}$
satisfying $|\Supp(\xi)|=m$.

Now, the intersection of $\Cal Z_{n, i+1}$ and $p_2^{-1}(\Gamma)$
contains an open dense subset $U_0$ which consists of all the
points $(\xi, x) \in \hil{n} \times X$ satisfying the conditions:
$$\ell(\xi_x) = i+1, \xi-\xi_x \in X^{[n-i-1]}_0, \, \Supp(\xi)
\cap \Gamma = \{x\} \text{ if } s > 0.$$ Furthermore, the
intersection $\Cal Z_{n, i+1} \cap p_2^{-1}(\Gamma)$ is
transversal along $U_0$. Since either $i > 0$ or $s > 0$, the
restriction $p_1|_{U_0}$ maps $U_0$ homeomorphically to the subset
$V_0$ of $\hil{n}$ consisting of all the points $\xi \in \hil{n}$
satisfying the conditions: $$\hbox{$\exists x\in \Gamma$ with
$\ell(\xi_x) = i+1$, $\xi-\xi_x \in X^{[n-i-1]}_0$}, \, \Supp(\xi)
\cap \Gamma = \{x\} \text{ if } s > 0.$$ It follows that the
Poincar\'e dual of $p_{1*}([\Cal Z_{n, i+1}] \cdot p_2^*\gamma)$
is the homology class represented by the closure of $V_0$, denoted
by ${\overline {V_0}}$: $${\overline {V_0}} = \{\, \xi\in
\hil{n}\, |\, \hbox{$\exists x\in \Gamma$ with $\ell(\xi_x)\ge
i+1$}\,\}.$$

Similarly, using (2.10) and induction, we conclude that the
Poincar\'e dual of $$B_i(\gamma, n) = {1 \over (n-i-1)!} \cdot
\frak q_{i+1}(\gamma) \frak q_{1}(1_X)^{n-i-1} |0\rangle$$ is also
represented by ${\overline {V_0}}$ (see \cite{Na2}). So Lemma 3.28
follows.
%$$R= \{\, \xi\in \hil{n}\, |\, \exists \eta\in \hil{n-i-1},
%\eta\subset \xi,\text{Supp}(I_{\eta}/I_{\xi})=\{x\}\in \Gamma\,\}.$$
\endproof

\proposition{3.29} Let $i \ge 0$, $n \in \Zee$, and $\gamma,
\alpha \in H^*(X)$. Then the commutator $[\frak B_i(\gamma), \frak
q_n(\alpha)]$ satisfies the transfer property, i.e., we have
$$[\frak B_i(\gamma), \frak q_n(\alpha)] = [\frak B_i(1_X), \frak
q_n(\gamma\alpha)] = [\frak B_i(\gamma\alpha), \frak q_n(1_X)].$$
\endstatement
\proof First of all, we notice that it suffices to show that
$$[\frak B_i(\gamma), \frak q_n(\alpha)] = [\frak B_i(1_X), \frak
q_n(\gamma\alpha)].   \tag 3.30$$ Next, recall that $\frak
q_0(\alpha) = 0$. Also, by (2.7), (2.11) and (3.4), we have
$$[\frak B_i(\gamma), \frak q_{-n}(\alpha)] = [\frak
B_i(\gamma)^\dagger, (-1)^n \cdot \frak q_n(\alpha)^\dagger] =
(-1)^{n+1} \cdot [\frak B_i(\gamma), \frak q_n(\alpha)]^\dagger.$$
So we need only to prove (3.30) for $n > 0$. In the following, let
$n > 0$.

Consider the action of $[\frak B_i(\gamma), \frak q_n(\alpha)]$ on
$a \in H^*(\hil{m})$. By the definition of the operators $\frak
B_i(\gamma)$ and $\frak q_n(\alpha)$, we see that $\frak
B_i(\gamma) \frak q_n(\alpha)(a)$ is equal to $$(p_{n+m,
1})_*([\Cal Z_{n+m, i+1}] \cdot (p_{n+m, 2})^*\gamma) \cdot \w
p_{1*}([Q^{[m+n,m]}] \cdot \w \rho^*\alpha \cdot \w p_2^*a)$$
where $p_{n+m, 1}, p_{n+m, 2}$ are the projections of $\hil{m+n}
\times X$ to $\hil{m+n}, X$ respectively, and $\w p_1, \w \rho, \w
p_2$ are the projections of $\hil{m+n} \times X \times \hil{m}$ to
$\hil{m+n}, X, \hil{m}$ respectively. Using the projection formula
and pulling all the cohomology classes $[\Cal Z_{n+m, i+1}],
(p_{n+m, 2})^*\gamma, [Q^{[m+n,m]}], \w \rho^*\alpha, \w p_2^*a$
to $X^{[n+m], [1], [1], [m]}$, we conclude that $$\frak
B_i(\gamma) \frak q_n(\alpha)(a) = p_{1*}(p_{12}^*[\Cal Z_{n+m,
i+1}] \cdot p_{134}^*[Q^{[m+n,m]}] \cdot p_2^*\gamma \cdot
p_3^*\alpha \cdot p_4^*a)  \tag 3.31$$ where for $1 \le i_1 <
\ldots < i_s \le 4$, the map $p_{i_1 \ldots i_s}$ stands for the
projection of $X^{[n+m], [1], [1], [m]}$ to the product of the
$i_1$-th, $\ldots$, $i_s$-th factors. Similarly, $$\frak
q_n(\alpha)\frak B_i(\gamma)(a) = p_{1*}(p_{24}^*[\Cal Z_{m, i+1}]
\cdot p_{134}^*[Q^{[m+n,m]}] \cdot p_2^*\gamma \cdot p_3^*\alpha
\cdot p_4^*a).  \tag 3.32$$ \noindent {\bf Claim.} {\it Regard
$p_{12}^*[\Cal Z_{n+m, i+1}] \cdot p_{134}^*[Q^{[m+n,m]}]$ and
$p_{24}^*[\Cal Z_{m, i+1}] \cdot p_{134}^*[Q^{[m+n,m]}]$ to be the
products of algebraic cycles in the Chow ring $A^*(X^{[n+m], [1],
[1], [m]})$. Then, $$\align &\quad (p_{12}^*[\Cal Z_{n+m, i+1}]
-p_{24}^*[\Cal Z_{m, i+1}])
 \cdot p_{134}^*[Q^{[m+n,m]}] \\
&\in j_*(A^*(X^{[n+m]}\times \Delta_X \times \hil{m})) \subset
 A^*(X^{[n+m], [1], [1], [m]})  \tag 3.33
\endalign$$
where $\Delta_X$ stands for the diagonal in $X^2 = X \times X$,
and $j$ is the inclusion $$X^{[n+m]}\times \Delta_X \times \hil{m}
\hookrightarrow X^{[n+m], [1], [1], [m]}.$$} \proof Denote a point
in $X^{[n+m], [1], [1], [m]}$ by $(\xi, x, y, \eta)$. Let $U
\subset X^{[n+m], [1], [1], [m]}$ be the open subset consisting of
all the points $(\xi, x, y, \eta)$ with $x \ne y$. Then the
complement of $U$ is precisely $X^{[n+m]}\times \Delta_X \times
\hil{m}$.

Consider $p_{12}^{-1}\Cal Z_{n+m, i+1} \cap
p_{134}^{-1}Q^{[m+n,m]}$ which has the expected dimension
$$2m+n+1-i.  \tag 3.34$$ A point $(\xi, x, y, \eta) \in U \cap
(p_{12}^{-1}\Cal Z_{n+m, i+1} \cap p_{134}^{-1}Q^{[m+n,m]})$ if
and only if $x \ne y$, $\xi_x = \eta_x$ has length greater than or
equal to $(i+1)$, and $I_\eta/I_\xi$ has length $n$ and support
$\{ y \}$, i.e., if and only if $$\xi = \eta_x +\xi_y + \zeta
\quad \text{and} \quad \eta = \eta_x +\eta_y + \zeta   \tag 3.35$$
where $x \ne y$, $\ell(\eta_x) \ge (i+1)$, $\eta_y \subset \xi_y$,
$\ell(\xi_y) = n + \ell(\eta_y)$, $\{ x, y\} \cap \Supp(\zeta) =
\emptyset$, and $\eta_x, \eta_y, \xi_y$ are supported at $x, y, y$
respectively. If $\eta_y \ne \emptyset$, then the dimension of the
set of those points $(\xi, x, y, \eta)$ satisfying (3.35) is at
most $$\align
 &\# (\hbox{moduli of }x, y) + (\ell(\eta_x)-1)+(\ell(\eta_y)-1)
  +(\ell(\xi_y)-1) +2\ell(\zeta) \\
=&2m+n+1-\ell(\eta_x) < 2m+n+1-i.  \tag 3.36
\endalign$$
If $\eta_y = \emptyset$, then the dimension of the set of the
points $(\xi, x, y, \eta)$ satisfying (3.35) is $$\align
 &\# (\hbox{moduli of }x, y) + \# (\hbox{moduli of } \eta_x)
  +(\ell(\xi_y)-1) + 2\ell(\zeta) \\
=&2m+n+1-\# (\hbox{moduli of } \eta_x) = 2m+n+1-i.  \tag 3.37
\endalign$$
By (3.34), (3.36) and (3.37), $U \cap (p_{12}^{-1}\Cal Z_{n+m,
i+1} \cap p_{134}^{-1}Q^{[m+n,m]})$ contains the open dense subset
$V \, {\overset \text{def} \to =} \, \{ (\xi, x, y, \eta) \, | \,
(\xi, x, y, \eta) \text{ satisfies (3.35) with } \eta_y =
\emptyset \}$ which is also irreducible. Now since $x \ne y$, the
intersection $p_{12}^{-1}\Cal Z_{n+m, i+1} \cap
p_{134}^{-1}Q^{[m+n,m]}$ along $V$ is transversal. So using the
refined intersection \cite{Ful}, we conclude that $$p_{12}^*[\Cal
Z_{n+m, i+1}] \cdot p_{134}^*[Q^{[m+n,m]}] = [{\overline V}] +
j_*(b_1)  \tag 3.38$$ where ${\overline V}$ is the closure of $V$
in $X^{[n+m], [1], [1], [m]}$ and $b_1 \in A^*(X^{[n+m]}\times
\Delta_X \times \hil{m})$.

Similarly, we see that $p_{24}^*[\Cal Z_{m, i+1}] \cdot
p_{134}^*[Q^{[m+n,m]}] = [{\overline V}] + j_*(b_2)$ for some
algebraic cycle $b_2 \in A^*(X^{[n+m]}\times \Delta_X \times
\hil{m})$. Combining this with (3.38) yields (3.33).
\endproof

Now we continue the proof of (3.30). Recalling our conventions
established at the end of section 1, we see from (3.31), (3.32)
and (3.33) that $$\align
 &[\frak B_i(\gamma), \frak q_n(\alpha)](a)  \\
=&p_{1*}((p_{12}^*[\Cal Z_{n+m, i+1}] -p_{24}^*[\Cal Z_{m, i+1}])
  \cdot p_{134}^*[Q^{[m+n,m]}] \cdot p_2^*\gamma \cdot
  p_3^*\alpha \cdot p_4^*a)   \\
=&{\overline p}_{1*}(b \cdot {\overline \rho}^*(\gamma \alpha)
  \cdot {\overline p}_2^*a)  \tag 3.39
\endalign$$
where $b \in H^*(X^{[n+m]}\times \Delta_X \times \hil{m})$ is
independent of $\gamma, \alpha$ and $a$, and ${\overline p}_1,
{\overline \rho}, {\overline p}_2$ are the projections of
$\hil{m+n} \times \Delta_X \times \hil{m}$ to $\hil{m+n}, \Delta_X
\cong X, \hil{m}$ respectively. Since $b$ is independent of
$\gamma, \alpha$ and $a$, (3.30) follows immediately from (3.39).
\endproof

\bigskip\noindent
{\bf Acknowledgments:} We thank the referee for his helpful
comments. Li's work partially supported by Hong Kong University of
Science and Technology grant number HKUST6170/99P. Qin's work
partially supported by National Science Foundation grant number
DMS-9996346 and an Alfred P. Sloan Research Fellowship. Wang's
work Partially supported by National Science Foundation grant
number DMS-0070422 and a Faculty Research and Professional
Development Fund at North Carolina State University.

\enddocument